\documentclass[11pt,letterpaper]{amsart}
\usepackage[utf8]{inputenc}
\usepackage{comment}
\usepackage{xcolor}
\usepackage[T1]{fontenc}

\usepackage{amsmath, amsfonts,amsthm,graphicx,fullpage,etex,tikz,hyperref,amssymb}

\let\arxiv\arXiv

\usepackage[capitalize]{cleveref}
\crefname{ineq}{Ineq.}{inequalities}
\creflabelformat{ineq}{#2{\upshape(#1)}#3} 
\usepackage[normalem]{ulem}
\usepackage{enumerate}
\usepackage{enumitem}
\usepackage{todonotes}
\usepackage{tcolorbox}
\hypersetup{
    colorlinks=true,
    linkcolor=blue,
    filecolor=magenta,      
    urlcolor=cyan,
}

\newtheorem{theorem}{Theorem}[section]

\newtheorem{Definition}[theorem]{Definition}

\theoremstyle{remark}
\newtheorem{remark}[theorem]{Remark}

\numberwithin{equation}{section}

\newcommand{\T}{\mathrm{T}}

\newcommand{\Z}{\mathbb{Z}}

\newcommand{\RP}{\mathbb{R}\mathbf{P}}

\newcommand{\R}{\mathbb{R}}

\newcommand {\Ad}{\mathrm{Ad}}
\renewcommand {\T}{\mathbb{T}}

\newcommand{\holder}{H\"{o}lder }
\newcommand{\smooth}{ C^{\infty}}

\def\rh{Rodriguez Hertz }

{\fontfamily{lmr}\selectfont}

\newtheorem{Structural Stability Theorem}[theorem]{Structural Stability Theorem}

\begin{document}

\title{On the work of Zhiren Wang on  Rigidity  in Dynamics}

\author{RALF SPATZIER$^\ast$}
\address{DEPARTMENT OF MATHEMATICS, UNIVERSITY OF MICHIGAN, ANN ARBOR, MI, 48109.}
\email{spatzier@umich.edu}

\thanks{$^{\ast }$ Supported in part by  NSF grants  DMS-2003712.}

\date{\today }

\date{}




\begin{abstract}
 In honor of Zhiren Wang on the occasion of being awarded the    Brin Prize, we report on his exciting and deep work on rigidity of higher rank abelian groups and lattices in higher rank semisimple groups.   
     \end{abstract}

\maketitle

\section{Introduction}   \label{section:Introduction}

\subsection{Panegyric} Zhiren Wang is the winner of the eleventh Brin Prize in Dynamical Systems.   He has made fabulous contributions to rigidity properties of actions by higher rank abelian and semisimple Lie groups and their lattices. This is complemented by remarkable work on Sarnak's disjointness conjecture.

Much of Wang's work deals with the classification of actions of a discrete group $\Gamma$ such that at least one element acts as an Anosov diffeomorphism, a so-called Anosov action of $\Gamma$. 
The main themes of this survey will be the following superb developments explicitly mentioned in the awarding of the Brin Prize:
\begin{enumerate}
    \item  smooth classification of Anosov actions of $\Z ^d, d \geq 2$, on tori and nilmanifolds (joint with F. Rodriguez Hertz) \cite{RH-W},
    \item topological rigidity of actions on nilmanifolds by higher rank lattices $\Gamma$ in semisimple Lie groups with hyperbolic linearization and  smooth classification of  actions of $\Gamma$ on nilmanifolds with an Anosov element (joint with A. Brown and F. Rodriguez Hertz) \cite{MR3702679},
    \item measurable classification and structure of higher rank lattice actions in low and intermediate dimensions  (with A. Brown and F. Rodriguez Hertz)~\cite{BFHW}.
\end{enumerate}

It is a great pleasure to report on these crowning contributions to the rigidity program and survey the landscape around them.  I will concentrate on Wang's work on actions of higher rank lattices in semisimple Lie groups, in part because I discussed global rigidity of Anosov $\Z ^d$ actions on tori and nilmanifolds in \cite{brinprize2015}.  

\subsection{The actors}
Let us begin by introducing our central players, certain discrete groups $\Gamma$ embedded into an ambient Lie group $G$.  Typically they are either $\Z ^k$, diagonalizable in $G$ (if $G \subset GL(N,\R)$ for some $N$), or $SL(k+1, \Z)$ for $k \geq 2$.  More generally, we can consider 
lattices $\Gamma$ in connected Lie groups $G$, i.e., discrete subgroups such that $G/\Gamma$ has finite Haar measure. To be precise,  we need to impose an additional condition, namely: we consider lattices $\Gamma$ in a higher rank semisimple Lie group with all factors of real rank at least 2, such as $SL(k+1, \R)$ for $k \geq 2$.  Indeed, $SL(k+1, \Z)$ is a lattice in $SL(k+1, \R)$ \cite{raghunathan1972}. 
Typically we assume that such $\Gamma$  are irreducible, i.e., that no finite index subgroups are non-trivial products.  Henceforth we will call such $\Gamma$ {\em higher rank lattices}  for short.

Most importantly, all of these groups contain $\Z^2$ as subgroups.  This follows from Dirichlet's unit theorem for $SL(k+1, \Z)$ for $k \geq 2$. For general higher rank lattices, this is a result of Prasad and Raghunathan \cite{Prasad-Raghunathan1972}.    Using the structure of the dynamics of these higher rank abelian subgroups will be the key to many results.

\subsection{Act I: Actions of large groups} 
R. Zimmer started investigating actions of higher-rank semisimple Lie groups and their lattices by diffeomorphisms on compact manifolds in the early 1980s. This comprises the {\em Zimmer Program}. 
He was much inspired by Margulis' superrigidity  theorems which classify finite-dimensional representations of higher rank lattices \cite{Margulis-ICM1975,margulis-book}.   Zimmer extended this work in his superrigidity theorem for cocycles around 1980 \cite{Zimmer-Annals1980,Zimmer1}.  He 
also made far reaching conjectures about such actions.  

One of these, the so-called {\em Zimmer conjecture},  suggested that lattices cannot act on low dimensional manifolds except via finite groups. This conjecture had been proven for one-dimensional manifolds and surfaces (under additional assumptions), cf., e.g., \cite{MR1259430,MR1911660,MR1946555,MR2219247,MR1911660,MR4521507}.
 In  recent breakthroughs, this has finally been proven in many cases---and in particular for $SL(n+1,\Z)$, $n \geq 2$---in three fabulous works by Brown, Fisher and Hurtado \cite{Brown-Fisher-Hurtado,MR4132960,Brown:2021wc}.   Note however that in general the range of the dimension of the underlying manifold is not yet optimal.
  Finding a $\Gamma$-invariant metric with good smoothness played a critical role.  Indeed Zimmer already realized that the superrigidity theorem for cocycles provides a measurable such metric, at least in the presence of a $\Gamma$-invariant volume.  Additional tools and insights however were needed to prove regularity of such metrics, thus getting an embedding into the isometry group of a compact Riemannian manifold which is compact itself. By Margulis  \cite{Margulis-book1991}, homomorphisms of $\Gamma$  into compact Lie groups are well understood.  Vice versa, such representations yield actions. Note that dynamically such actions are not very interesting exactly as they are isometric.

\smallskip

On the other end of the dynamical spectrum lie uniformly hyperbolic actions.   Here the dynamics gets complicated, especially the structure of the orbits.  Closed orbits coexist with dense ones. Lyapunov exponents are typically positive. One has mixing and even exponential mixing. In the uniformly hyperbolic case,  dynamics typically comes  with \holder structures. One can then  hope to   use \holder regularity when trivializing the derivative cocycle, as Zimmer already realized. Again, it took many years and advances in other areas, smooth dynamics in particular, to overcome the difficulties.  As we will see, Wang's work with A.~Brown and F. Rodriguez Hertz on such actions on tori and nilmanifolds 
plays a central role \cite{MR3702679}.  Indeed, their work gives an essentially complete classification of smooth actions with an Anosov element on nilmanifolds.  Miraculously, they even obtain $C^0$ results essentially when the ``linearization'' of an action of $\Gamma$ contains an  element without eigenvalues of modulus 1.  In the case of an action on a torus, this linearization is simply given by the induced action on first homology.

\smallskip

Finally, there are also projective actions, in particular the natural, classical ones for higher rank lattices in $SL(n,\R)$ on real projective space, Grassmannians and other flag varieties.  More generally one gets  actions on boundaries $G/Q$ where $Q$ is a parabolic in a semisimple group $G$.  One can easily build  actions of higher rank lattices and semisimple groups by a suspension construction from actions of a parabolic subgroup. As any flow will give rise to an action of a parabolic subgroup, such actions cannot   cannot be classified in a real sense.   Note also that these actions do  not preserve any probability measure. Naturally the question arises if any action is built from projective actions in some fashion, or at least factor over them.  Nevo and Zimmer have two papers \cite{NevoZ00,NevoZ02} with measurable versions of such statements.  Unfortunately, their conditions are difficult to use.  Here again, Wang with Brown and Rodriguez Hertz has much clarified the picture in      \cite{BFHW}.  For lattices $\Gamma$ in   $SL(n,\R)$ and  $n \geq 3$ say, they can classify  smooth actions of such lattices when the dimension of the manifold is $n-1$.    If the action preserves a volume, they can even handle dimension $n$.  
Furthermore they establish  a factorization picture for actions in intermediate dimensions between $n$ and $2n$.  Currently, they are  developing a smooth version of these results.

\smallskip

Dynamics is often said to have three main characters, elliptic, hyperbolic and parabolic dynamics.  And so it is in the Zimmer program.  Really, these situations are the extremes and much work remains to be done.  

\subsection{Act II: Hyperbolic actions of higher rank abelian groups}
One of the key engines in the Zimmer program is the analysis of higher rank abelian groups with special dynamical properties. In fact, the Zimmer program provided much of the motivation for studying the latter. It holds its own interest as a study in symmetry in dynamics.   Again one strives to find strong structural properties, with the eventual goal of classifying  such actions. Best understood is the uniformly hyperbolic case for $\Z ^k$ actions.  If the action is on a torus or nilmanifold $N/\Lambda$, Wang and Rodriguez Hertz show in \cite{RH-W}  that the action is smoothly conjugate to an action by automorphisms supposing that its linearization does not factor (algebraically) over a $\Z$ action.  Their theorem   resolves a well-known conjecture of Katok and this author for uniformly hyperbolic actions of higher rank abelian groups if the underlying manifold is a nilmanifold. As the hypotheses only involve knowledge of the underlying manifold, we can interpret this as a {\em global rigidity} statement.  This is in contrast to {\em local rigidity} where one considers actions $C^1$-close to a given one on  a fixed finite generating set. 
All of these advances are part of a larger program to analyze and potentially classify actions of higher rank abelian groups assuming sufficient hyperbolicity for the action.

We will now explore these topics in greater detail,  starting with Wang's paper on $\Z^k$  actions.   We refer to  \cite{brinprize2015} for an extensive discussion of related conjectures as well as some history and motivation.

\section{Anosov actions by abelian groups} 

\subsection{Some background}

Let $\Gamma$ be a discrete group, and suppose 
 $\Gamma$ acts on a compact manifold $M$ by $C^1$-diffeomorphisms.  We call the action {\em Anosov} if some element $\gamma _0 \in \Gamma$ is an 
Anosov diffeomorphism of $M$.  We recall the definition:

\begin{Definition}
Let $M$ be a compact manifold with some ambient Riemannian metric $\|\:\|$. Then a diffeomorphism  $f: M \to M$ is called {\em Anosov} if there is a continuous splitting of the tangent bundle $TM = E^s \oplus E^u$  and constants $C>0$ and $\lambda >0$ such that for all $n >0$, $v \in E^s$  and $w \in E^u$
\[  \|df ^n (v)  \|  \leq C e^{-n \lambda }  \|v\|  \:\:\:\text{ and }  \:\:\: \|df ^{-n} (w)  \|  \leq C e^{-n \lambda }  \|w\|   .\]
\end{Definition}

Once you have an Anosov diffeomorphism  $f: M \to M$, there are many, thanks to structural stability \cite{KatokHasselBook}.  Indeed, if a diffeomorphism $g: M \to M$ is sufficiently close to $f$ in the $C^1$-topology, then $g$ is Anosov.  Moreover, there is a \holder homeomorphism $\phi$  homotopic and $C^0$ close  to $id _M$ which conjugates  $g$  to $f$,
 $  g = \phi \circ f \circ \phi ^{-1}$.  We can think of structural stability as a rudimentary rigidity theorem.  Note that the \holder conjugacy $\phi$ is no better than \holder in general.  In particular, it is not $C^1$ typically as is easy to see by perturbing the map near a periodic point, changing the derivative. 
Examples of Anosov diffeomorphisms arise from linear maps on tori without eigenvalues of modulus 1, and more generally automorphisms of nilmanifolds $N/\Lambda$ where $N$ is a simply connected nilpotent group and $\Lambda \subset N$ a cocompact lattice. It is an open problem which nilmanifolds admit Anosov automorphisms.  For example the Heisenberg group $H$ does not admit such quotients while $H \times H$ does, via the Borel-Smale examples \cite{Smale}. This problem has been studied by many authors, e.g., Dani, Lauret, Mainkar, Payne, Der\'{e}, Will  (e.g., \cite{MR518242,MR2140439,MR2481335,MR2471923,Dere:2020aa}).  They were classified in low dimensions.
More generally, compositions of translations and automorphisms on a nilpotent group $N$ may descend to a nilmanifold.  In that case, one calls the resulting map an  {\em affine } Anosov diffeomorphism.

Anosov and Smale asked already in 1967 whether any Anosov diffeomorphism $f$ is topologically conjugate to an affine Anosov diffeomorphism, possibly after passing to a finite cover of $f$
 \cite{Smale}.  In general, this is still a wide open problem.
Franks and later Manning made major advances to this problem in \cite{Franks1970, Manning74}.  In particular, they showed that the conjecture holds for Anosov automorphisms of infra-nilmanifolds.   The resulting conjugacy is often called the {\em Franks--Manning conjugacy} and the automorphism its {\em linearization}.   It is unique in the homotopy class of the identity.  In general, this conjugacy is only H\"{o}lder.   We describe the conjugacy in more detail below in \cref{sec-outline}.

\subsection{Anosov diffeomorphisms with nontrivial centralizers} The Franks--Manning conjugacy will automatically conjugate commuting diffeomorphisms as well, as is easy to see.  To obtain a smooth classification of Anosov actions for $\Z^k, k \geq 2$ now becomes a regularity problem:  Does the presence of a ``non-trivial'' commuting map make this conjugacy smooth? \rh~and Wang have the ultimate answer in \cite{RH-W}.  First, though, recall the notion of factor action of a $\Z^k$ action on a space $X$.  This is simply an action of $\Z^k$ on a space $Y$ and a $\Z^k$-equivariant map $\pi: X \to Y$.  We call a factor action {\em rank 1} if some $\Z ^{k-1}$ acts trivially on the factor. We can talk about factors in various categories, in particular measurable, continuous, $C^1$ or $C^{\infty}$, or even algebraic ones.

\begin{theorem} [\rh--Wang]			 \label{theorem:global rigidity}
Suppose $k \geq 2$ and that $\alpha$ is a   $\smooth$ $\Z ^k$ action on an infra-nilmanifold M with at least one Anosov element.   Suppose the linearization $\rho$ of $\alpha$ does not have (algebraic) rank 1 factors.  Then the Franks--Manning conjugacy is a $\smooth$ diffeomorphism.  
\end{theorem}  

This result is particularly pleasing as existence of an algebraic rank one factor of the linearization is a checkable condition.  It is also the culmination of a long series of works on local and global rigidity of uniformly hyperbolic actions of $\Z^k, k \geq 2$ on tori and nilmanifolds,  These started with Hurder \cite{Hurder,hurder-survey1994} who proved various local and deformation rigidity rigidity results for $\Z ^k$ Anosov actions on tori and nilmanifolds by automorphisms if the maximal non-trivial intersections of stable manifolds of different Anosov elements of $\Z ^k$ are one-dimensional.  Here local rigidity of an action $\alpha$ means that if another action $\beta$ is sufficiently $C^1$-close to $\alpha$ on some fixed finite generating set of $\Z ^k$, then $\beta$ is smoothly conjugate to $\alpha$.  Deformation rigidity is weaker and  assumes additionally that $\beta$ is connected to $\alpha$ by a continuous  path of smooth actions. From his work on $\Z^k$ actions, Hurder deduced deformation rigidity of the standard action of $SL(n,\Z)$ on the $n$-torus \cite{Hurder}.  This was quickly followed  by works of Katok, Lewis and in part Zimmer for local rigidity for more general actions \cite{KL1, KL2} and even global rigidity under special assumptions \cite{KL}.  Local rigidity of certain algebraic higher rank Anosov actions of $\Z ^k$ and also $\R ^k$ was proved by Katok and this author in \cite{KS97}.

The immediate precursor to the work of \rh~and Wang was achieved  by 
D. Fisher, B. Kalinin and this author in  \cite{FKS} under the additional assumption   that $\alpha$ has many Anosov diffeomorphisms.   As we will see, both this result and the methods used will play an important role in the proof of Theorem \ref{theorem:global rigidity}.  The basic strategy of \rh~and Wang is to show that there are lots of Anosov elements, thus reducing their result to the earlier one in \cite{FKS}.   This may sound easy, but is actually hard and requires  {\em  great finesse}. 

\subsection{Outline of the arguments} \label{sec-outline}
Let us first summarize the arguments from \cite{FKS}.  The main point is that  the Franks--Manning conjugacy $h$ is  explicitly given as a  series in terms of projections to eigenspaces of the linearized action by: 
\begin{equation} \label{hVseries}
h_V (x) = \sum _{m=0} ^{\infty} A_V^{-m }\:Q_V (a^{m} x).
\end{equation}

This series always converges in the space of \holder functions, even for a single Anosov diffeomorphism---that is Franks' insight.  To see that $h$  is smooth in the higher rank setting of a $\Z ^k$ action, one studies its regularity along  the so-called {\em coarse Lyapunov spaces} of the unknown action. These are  maximal non-trivial intersections of stable manifolds of Anosov elements.  For the linearized action, Lyapunov exponents are simply logarithms of absolute values of eigenvalues. 
One then treats the Lyapunov exponents   as linear functionals on $\R ^k$.  The coarse Lyapunov spaces  become  sums of the Oseledets spaces of positively proportional Lyapunov functionals.  One can do a similar construction for an arbitrary $\Z ^k$ action.  However, the Oseledets spaces will only depend measurably on the base point.  To get \holder dependence, one needs to assume that there are enough Anosov elements in the action so that the measurable coarse Lyapunov spaces become intersections of stable manifolds.  For future reference, we will call the kernels of the Lyapunov functionals {\em Weyl chamber walls}.  Then {\em Weyl chambers} are  the connected components of the complement of the union of the Weyl chamber walls, just like in Lie theory. 

 Writing the expression for the first derivative of the Franks--Manning conjugacy, the contributions from the derivative of $a^m$ and the linear map $A_V^{-m } $ in this series balance each other out, and we cannot control convergence.   However, a more careful analysis shows that this series will still converge as a distribution, dual to \holder functions. Actually, one considers a non-Anosov element $\bar{a}$ and corresponding $\bar{A}$ in the suspension $\R^k$ action for which $a$ lies on a Weyl chamber wall, i.e., $\bar{A}$ has eigenvalues of modulus 1.  Then one does not have convergence of the series from traditional dynamics.  Rather one uses exponential mixing of the $\Z ^k$ action on the nilmanifold as established for toral automorphisms in \cite{FKS} and for nilmanifolds in \cite{MR3413978}.   In fact, partial derivatives of any order of the distribution along coarse Lyapunov foliations will converge. Using PDE techniques, one can then prove that $h$ is $\smooth$.  This requires the coarse Lyapunov foliations to vary \holder transversely and have smooth leaves. By the above, this is the case when there are enough Anosov elements. 
 
 This is precisely where Wang and Rodriguez Hertz' come in: they actually prove that there are always enough Anosov elements if the linearized action does not have a rank 1 factor.  In terms of Weyl chambers, they show that each Weyl chamber has an Anosov element.  The Weyl chamber walls clearly are not Anosov as some  Lyapunov export is 0 on each wall.  Starting with a Weyl chamber wall  $\mathcal{C}$ with an Anosov element, they show that any adjacent Weyl chamber wall also has an Anosov element.  They use   Ma\~{n}e's    characterization of Anosov diffeomorphisms  \cite{Mane77}: the dimensions of the stable manifolds at all periodic points need to be  the same and for all non-zero tangent vectors $v$ one has 
$\sup _{i \in \Z}\|d f^i (v) \| = \infty $.
Since the Franks--Manning conjugacy will map stable manifolds to stable manifolds, the dimension statement follows easily.  
Also, the growth statement is easy to prove at points at which Oseledets spaces exist for the action for the pull-back of Haar measure under the conjugacy. Proving the second claim is several magnitudes more difficult, and requires existence and regularity properties of the Franks--Manning conjugacy. 
Here is the short version of their argument: they show regularity of the Franks--Manning conjugacy, more precisely of the projection $h^0 _a$ to the slow stable on the algebraic side for some suitable element $a$, again using exponential mixing.  Then they show that $h^0 _a$ has constant rank along all stable manifolds.  This is difficult, and technical.  To give a flavor, they suppose the singular set is non-empty, and then analyze an invariant measure on the singular set comparing it with the pull-back of Lebesgue measure by the Franks--Manning conjugacy. We refer to \cite{RH-W} or \cite{brinprize2015} for more details. 

\subsection{Comments, problems and related results} Let us briefly comment on later developments and possible conjectures.

\medskip

{\bf A:} If the Anosov-Smale conjecture proves true, then we automatically get a classification up to smooth conjugacy of  $\Z ^k$ actions on any closed manifold $M$.  Indeed, $M$ then is homeomorphic to an infra-nilmanifold by a conjugcay $\varphi$. Using differential topology \cite[Appendix by Jim Davis]{FKS}, a finite cover of $M$ is diffeomorphic to $N/\Lambda$.  One can then find a possibly larger finite cover to which the $\Z ^k$ action lifts.  Applying \cite{RH-W} gives a smooth conjugacy to an affine algebraic action assuming that the original action does not have continuous rank 1 factors.   Finally $\varphi$ inherits the smoothness from the lift. 

\medskip

{\bf B:}  Lei Yang and this author classified actions on closed manifolds by commuting expanding maps up to smooth conjugacy, again if there are no continuous rank one factors \cite{Spatzier-Yang2017}.   For expanding maps the analog of the Anosov-Smale conjecture is known thanks to work by Gromov and Shub \cite{GromovPolyGrowth,ShubExpandingMaps1970}.    Then one can proceed in a similar way to the Anosov case after establishing the needed exponential mixing for the solenoid extension (for the latter, $p$-adic arguments are needed).  

\medskip

{\bf C:} Based on Wang's and \rh's work, one might  speculate that existence of one Anosov element will always force  existence of lots, even  dense set of Anosov elements.  For $\R^k$ actions however, K. Vinhage recently constructed examples of Anosov $\R^k$ actions for which the Anosov elements do not form a dense subset \cite{Vinhage:2022}.  His construction involves a time change of a product $\R^k$ action. Potentially, all such examples may arise from time changes.  
For $\Z^k$ actions, having one Anosov element may prove sufficient. 

\medskip

{\bf D:}. For $\R^l \times \Z ^d$ actions, $d + l \geq 2$, on arbitrary closed manifolds with a dense set of Anosov elements and transitive one parameter subgroups, Vinhage and the author  proved a smooth classification theorem under the additional assumption that  intersections of stable foliations of suitable elements are 1-dimensional.  We call such actions {\em Cartan} actions \cite{Spatzier:2020tx}.    Part of this classification has been extended to general higher rank Anosov actions \cite{Damjanovic:2022aa}.  This follows earlier work, under much more stringent assumptions, by Damjanovic and Xu as well as Kalinin and Sadovskaya and this author \cite{MR4145817, MR2240907,KalSad07,KaSp04}. 

\medskip

{\bf E:} These works give support to a conjecture by Katok and this author, essentially asserting  smooth classification of higher rank $\Z^k$ or $\R^k$ actions with a dense set of Anosov elements.

\smallskip

We note that the case of $\R^k$ actions is significantly different from $\Z ^k$ actions since  there is no reasonable conjectural classification of Anosov flows. Indeed, given a closed manifold $N$ of negative sectional curvature, its geodesic flow on the unit tangent bundle $SN$ is Anosov.  Many such $N$ have been constructed by Ontaneda in \cite{Ontaneda:2020aa}.   Moreover, there are now many non-geometric examples of Anosov flows, on many different underlying manifolds \cite{Barthelme:2021aa,Bonatti:2023aa}.
\vspace{.5em}

{\bf F:} While still in its  infancy, much is already known for partially hyperbolic diffeomorphisms and their  centralizers, generalizing the Anosov case.  Naturally, the situation  is much more complicated, as explored by Damjanovic, Wilkinson and Xu \cite{Damjanovic:2019aa,Damjanovic:2023aa}.

\section{Passage from  discrete groups to Lie groups}

Let us first introduce some very general ideas and constructions. They  will play a key role in the discussion of actions of lattices in Sections
\ref{sec:anosovlattice} and \ref{sec:latticelowdim}.

Actions by discrete groups are typically difficult to analyze.  Actions by  Lie groups on the other hand are considerably simpler.  For example, one can study subgroups  satisfying particular dynamical properties such as subexponential growth in certain directions.  For the discrete groups this may only be possible approximately,  along suitable sequences which is more difficult then along a Lie subgroup. 

 When $\Gamma$ is a lattice in a Lie group $G$, we can transition from an action by $\Gamma$ to one by $G$.  The latter is  called the {\em suspension} or also {\em induced action}.  In the simplest case of $\Z$ in $\R$, this is nothing but the suspension flow or mapping torus construction.  

 \begin{remark} In the following, the reader can always supplant a general semisimple Lie group $G$ by $SL(n,\R)$.  Then the {\em split Cartan subgroup} $A$ simply becomes the diagonal subgroup of $SL(n,\R)$.  

 We will consider linear representations $\rho$ of $G$ and study their {\em weights}.  To the dynamicist, these are nothing but Lyapunov exponents of the representation, thought of as linear functionals on $A$. 
 
 For the most basic representation, the adjoint representation by the derivatives of conjugation by $g \in G$ on the Lie algebra of $G$, we call the weights {\em roots}. The roots are basic for the structure theory of $G$ and its representations.  In the case of $G = SL(n,\R)$, the roots are nothing but the linear functionals $e_i - e_j$ where the $e_i$ are the $i$-th diagonal entry of a diagonal matrix. 

 Finally, the {\em root spaces} are just the eigenspaces of the roots, Oseledets spaces in dynamical terms.  For $G=SL(n,\R)$, these are just the $\R \: E_{ij}$ subspaces of the Lie algebra of trace 0 $n \times n$ matrices of $SL(n,\R)$ where $E_{ij}$ is the $i j$ elementary matrix with a 1 in the $i j $ position and 0 elsewhere.
 \end{remark}

 \subsection{Suspensions} \label{sec:susp}  Crucial to the authors' approach is the use of the suspension of the $\Gamma$ action to a $G$ action for the ambient semisimple Lie group $G$.  More precisely, given an action $\tau$ of $\Gamma$ on a space $X$, form the product $G \times X$, and mod out by the $\Gamma$ action defined by: $\gamma (g_0, x) = (g_0 \gamma ^{-1}, \tau (\gamma) x)$. Then $G$ acts on the resulting space $G \times _{\Gamma} X$ by left translation on the $G$-factor.  This gives a $G$ action $\hat{\tau}$ which is called the {\em suspension} (or also {\em induced action}) of the $\Gamma$ action.   This is the natural generalization of the suspension of a diffeomorphism to a flow. 

This construction allows  to make arguments for  $G$ actions which is decidedly easier than for  $\Gamma $ actions. In particular, one can study actions of subgroups of the split Cartan (i.e., diagonal ) subgroup. Of pertinent interest are the kernels  of linear functionals defined by the eigenvalues of the linearized action.  The derivatives of such maps grow only subexponentially in the eigenspace directions which will prove important later. In $\Gamma$ however, this subgroup may easily be trivial.

This simple construction already proved very useful when considering $\Z ^k$ Anosov actions, say on an $n$-torus $\T^n$.  Again, one would like to use elements in $\Z ^k$ which act isometrically in certain eigenspace directions - impossible in general for $\Z^k$ but trivial for the suspension.  In this case, the properties of the $\Z^k$ and $\R ^k$ actions are typically closely aligned.  For suspensions of actions of general $\Gamma$, this is much more subtle.  For example, it is not at all clear if the suspension of an Anosov $\Gamma$ action yields an  Anosov $G$ action, even when considering a uniform lattice $\Gamma \subset G$. 

 \subsection{Non-resonance} It will often be important to distinguish the base directions $G/\Gamma$ in the suspension from the fiber directions, in a dynamical fashion.  When we induce linear actions on tori or actions $\rho$ by automorphisms on nilmanifolds more generally, this is indeed possible.  
The crucial insight is that the representations $\rho$ and the adjoint representation $\Ad$ (on the Lie algebra by derivatives of conjugations) could be very different.  Indeed, call (restricted) non-zero weights $\alpha$ and $\beta$ of two representations $\rho _1$ and $\rho _2$ of $G$ {\em resonant} if $\alpha$ and $\beta$ are positively proportional as linear functionals on the Lie algebra of the split Cartan. Moreover, we call $\alpha $ and $\rho _2$ {\em non-resonant} if $\alpha$ is non-resonant with every non-zero weight of $\rho _2$. 

With regards to the adjoint representation, 
call a representation $\rho$ of $G\,${\em strongly non-resonant} if all non-zero (restricted) roots are non-resonant with  $\rho$.   Finally, the authors call $\rho$ {\em weakly non resonant} if the root spaces coming from the non-resonant (restricted) roots generate the Lie algebra of $G$ as a Lie algebra.   They then show that $\rho$ is weakly non-resonant provided that no (restricted) weight of $\rho$ is 0.  This is strictly an argument about finite dimensional representations of~$G$.

\smallskip

{\em Key feature of non-resonance}:  
Suppose $\Gamma$ acts on $X$, say a torus or nilmanifold for now.  Then for the suspension action,   we can ``freeze'' movement along a weight direction $\alpha$ in $X$ by acting by an element in $\ker \alpha$  and still ``move'' in all non-resonant root directions, i.e., contract or expand in those.

\smallskip

As will become clear below, this yields a connection  with the Anosov property:  Indeed, if a linear action (or action by automorphisms) $\rho$ is Anosov, then  no (restricted) weight of $\rho$ can be 0.  
Thus for an Anosov action, the toral (fiber) directions behave differently dynamically from the base directions for the suspension action. In other words, $\rho$ is strongly non-resonant.

\section{Anosov actions of higher rank lattices} \label{sec:anosovlattice}

We will call any action $\tau$ of a discrete group $\Gamma$ by $C^1$ diffeomorphisms {\em Anosov} if $\tau (\gamma)$ is an Anosov diffeomorphism for some element $\gamma \in \Gamma$.  The work on Anosov actions of higher rank abelian groups by Rodriguez Hertz and Wang strongly suggested that actions of higher rank lattices  $\Gamma$  on tori and nilmanifolds with an Anosov element should be classifiable.  Indeed, such  $\Gamma$  contain $\Z ^k$, $k \geq 2$, subgroups with Anosov elements.  By the abelian result, the latter will be smoothly conjugate to an algebraic action (if sufficiently irreducible).  However, it is not clear a priori how to manufacture a conjugacy that works jointly for all  of $\Gamma$ from the individual conjugacies for the higher rank abelian subgroups.  In the past, this was achieved in various ways, the most general appealing to Zimmer's cocycle superrigidity theorem, cf., e.g., \cite{KL2}.  However, this required existence of an invariant measure of full support for $\Gamma$ which again is not at all evident.  

In \cite{MR3702679}, Brown, Rodriguez Hertz and Wang overcome these problems and arrive at rather general and indeed best  results known classifying Anosov actions of higher rank lattices  on tori and nilmanifolds.  In the course of this work they also achieve a major breakthrough in the $C^0$ Zimmer program, i.e., the study of actions of higher rank  lattices on manifolds by homeomorphisms.  For such assuming the existence of an Anosov diffeomorphism does not make sense.  However, for tori, one can consider the induced action on first homology, arriving a a linear representation of $\Gamma$ called the {\em linear data} of the action of $\Gamma$.  Similarly, for nilmanifolds $M = N/\Lambda$ ($N$ a nilpotent simply connected group and $\Lambda$ a lattice in $N$), one gets an action by automorphisms on the (nilpotent) fundamental group $\pi_1 (M)$ which in turn gives a representation into the automorphism group of the ambient nilpotent group $N$ - unfortunately only up to conjugacy.   Here the technical assumption on lifting the action enters as it allows to  get a well defined representation into the automorphism group of $\Lambda$. Again they call this the {\em linear data} $\rho$ of the $\Gamma$ action.  They also require that the $\Gamma$ actions lift to the universal cover of $M$ which is often automatic, e..g for cocompact $\Gamma$.  We will not discuss this in detail and refer to \cite{MR3702679} instead.   Finally, if $\Gamma$ acts on two compact topological spaces $X$ and $Y$, we call a continuous surjective map $h: X \mapsto Y$ a {\em semi-conjugacy} if for all $\gamma \in \Gamma$ and all $ x \in X$, $h (\gamma x) = \gamma (h(x))$.  We call $h$ a {\em conjugacy} if $h$ is a homeomorphism. 

We are now ready to state the major breakthrough in \cite{MR3702679}.

\begin{theorem}[Brown--\rh--Wang]    \label{thm:semiconjugacy}
Suppose $\Gamma$ is a higher rank lattice acting $C^0$ on a  nilmanifold $M=N/\Lambda$ such that the linear data $\rho (a)$ for some $a \in \Gamma$ does not have eigenvalues of modulus 1.  Suppose further that the action lifts to the universal cover.   Then on a subgroup of finite index $\Gamma ' \subset \Gamma$, the $\Gamma '$ action is $C^0$ semi-conjugate to its linear data $\rho$ by a surjective continuous map $h:M \mapsto M$, homotopic to the identity.  
\end{theorem}

From this, they can solve the global rigidity problem for smooth actions of higher rank lattices on tori or nilmanifolds.

\begin{theorem}[Brown--\rh--Wang]    \label{thm:smoothrigidity}
Suppose $\tau$ is an Anosov action by a higher rank lattice  $\Gamma$ by $C^{\infty}$ diffeomorphisms on a  nilmanifold $M=N/\Lambda$. Suppose further that $\tau$ lifts to the universal cover.   Then the $\tau $ action is conjugate to its linear data $\rho$ by a  $C^{\infty}$ diffeomorphism  $h:M \mapsto M$, homotopic to the identity.  
\end{theorem} 

As a simple corollary, $C^1$ local rigidity of algebraic Anosov actions on nilmanifolds ensues since any $C^1$-close perturbation will still have an Anosov element by structural stability.

\subsection{Outline of the arguments} Theorem \ref{thm:smoothrigidity}  follows relatively easily from Theorems \ref{thm:semiconjugacy} and \ref{theorem:global rigidity}.
We will concentrate on the proof of Theorem \ref{thm:semiconjugacy}, and start with a review of basic constructions and definitions.

\subsection{Linear data}  Given any action of a group by homeomorphisms on a torus or even nilmanifold, one can associate a linear representation to it, called its {\em linear data}.  In the toral case, this is simply the induced action on the first homology group.  

Next, consider a nilmanifold $N/\Lambda$ where $N$  is a simply connected nilpotent group and $\Lambda$ is a lattice in $N$ (which is automatically cocompact).  These are well understood by work going back to Malcev \cite{Malcev:1951aa}, and can in fact be nicely described using the Lie algebra of $N$ \cite{MR3702679, GreenTao}.

Suppose $\Gamma$ acts by $\tau$ on $N/\Lambda$. Note that   $\Gamma$ naturally acts by outer automorphisms of $\Lambda$.  This action induces  a representation on $\Lambda $ and then $N$ (by Zariski density of $\Lambda $ in $N$) 
if we can pick a unique representative of the outer automorphism by an automorphism.  In particular, this is the case when the given $\Gamma$ action lifts to $N$.  We will always assume liftability,  similar to \cite{MR3702679}.  This representation of $\Gamma$ is called the  {\em linear data} of $\tau$.

\subsection{Continuous  semi-conjugacies} \label{sec:semi} We will now describe  the strategy of the  $C^0$ global rigidity result, Theorem \ref{thm:semiconjugacy} that the given action is $C^0$ semi-conjugate to its linear data.  The case of non-uniform lattices is much more technical since the non-compactness of the fundamental domain causes a variety of issues.  For simplicity for this exposition, we will assume that $\Gamma$ is co-compact.

\subsubsection{Toral case} Let us assume for starters that the underlying manifold of the $\Gamma$ action is a torus $\T^n$. 
As already mentioned above,   Franks proved the fundamental theorem in  \cite{Franks1970} from 1970 that any homeomorphism $f:  \T^n \mapsto \T^n$ is semi-conjugate to its linearization $A \in GL(n,\Z)$ provided that $A$ has a no eigenvalues on the unit circle, i.e., if $A$ is is Anosov (also cf. \cite[p.588]{KatokHasselBook}).   Thus we get a semi-conjugacy for $\rho(a)$ where $a \in \Gamma$  has  Anosov linear data. We described this in Equation \ref{hVseries}.

It is not hard to show that this semi-conjugacy extends to the centralizer of $a$ in $\Gamma$.  However, extending the semi-conjugacy to $\Gamma$ requires new ideas and is in fact one of the key achievements of this article.  First the authors work with the suspension actions of $G$ on $G \times _{\tau(\Gamma)} \T^n$ and $G \times _{\rho (\Gamma)} \T^n$ instead.  For these, they construct a (measurable) semi-conjugacy between these actions, as we will see just below for the general case of actions on nilmanifolds.


  
\subsubsection{Nilmanifold case} Again  consider a nilmanifold $N/\Lambda$ as above. As a first problem in the nilmanifold case, there is no precise analog for Franks' theorem about semi-conjugacies.   Manning did prove existence of a conjugacy but only when $\tau (a)$ is Anosov in \cite{Manning74}.   His arguments however do not apply to homeomorphisms and $C^0$ actions.  Instead, Brown, \rh~and Wang construct semi-conjugacies directly.  The idea is similar to Franks for tori, by writing a formal expression for the desired semi-conjugacy using a power series involving powers of $\rho(a)$, again as in Equation \ref{hVseries} via projecting to eigenspaces of $\rho (a)$ on the Lie algebra level.  If $\rho (a)$ is hyperbolic, we get convergence for the projection to the stable manifolds of $\rho (a)$, and also the unstable using $\rho (a) ^{-1}$. 

 Actually they want to construct a semi-conjugacy $\tilde{\Phi}$ between the suspension $G$ actions  of $\tau$ and $\rho$ on  $G \times _{\tau(\Gamma)} N/\Lambda$, first for $a$, then the centralizer of $a$ in $G$, and finally for $G$. Recall that the underlying space for the suspension is a fiber bundle with fibers $N/\Lambda$ over the base space $G/\Gamma$.  The desired semi-conjugacy $\tilde{\Phi}$ can then be viewed as a family of maps from $\tilde{\Phi} _g: N/\Lambda \mapsto N/\Lambda$  for $g \in G$  which intertwines correctly with the $\hat{\tau} (a)$ and $\hat{\rho} (a)$ suspension actions on the one hand and the $\Gamma$ action on the other. 

The basic idea is the same as for Franks' theorem: write down a formal solution for $a$, and show it converges.  In the end, they will get such a family of continuous maps depending measurably on $g \in G$ (the issue with measurability coming from a gluing construction involving fundamental domains).

As the target group is only nilpotent, not abelian, they actually do this by an inductive argument via the central series.  
More precisely, suppose the nilpotent group $N$ has nilpotency degreee $r$, and  let $N \mapsto N_1\mapsto \dots \mapsto N_r =\{1\} $ be a series of central extensions.  Then $\Lambda _i = N_i \cap \Lambda$ are lattices in $N_i$.  Thus 
one can think of the nilmanifold as a sequence of toral extensions starting with the maximal toral quotient $\T^n = (N/N_{r-1}) / (\Lambda/\Lambda _{r-1} )$.   For the latter,  one can write down formal expressions similar to Franks' case and prove convergence.   Thus they get a family of semi-conjugacies $\hat{\Phi} _1$ from $\tau (a)$ to the automorphism of $\T^n$ induced by $\rho (a)$.
One interprets $\hat{\Phi}$  as a map from $N/\Lambda $ to itself that is a semi-conjugacy for $a$ up to an error term in $N_{r-1}$. Then they improve the semi-conjugacy by an inductive argument on the degree of nilpotency subsequently converting the error terms on the various fibers into semiconjugacies. 
This yields the desired semi-conjugacy for $a$ on the suspensions.

\subsection{Extending the semiconjugacy  and the no resonance condition}  \label{sec:semiconj-nr} So far $\hat{\Phi}$ is a  semiconjugacy between $\hat{\tau} (a)$ to $\hat{\rho} (a)$.  It follows fairly easily that $\hat{\Phi}$ also provides a semi-conjugacy for the whole centralizer of $a$ in $G$.  It is far from trivial however that $\hat{\Phi}$ is a semi-conjugacy for $G$.  
Here the non-resonance between the roots of $G$ and the weights of $\rho$ enters crucially. More precisely, they need that the non-resonant roots generate the Lie algebra of $G$ as a Lie algebra. 

Consider a root vector $X$ in the Lie algebra of $G$ for the non-resonant root $\chi$.  The idea is simple.  If $\hat{\Phi}$ was differentiable, we could consider the pushforward vector field  $D \hat{\Phi} (X)$. If we move $X$ by $\ker \chi$, neither $X$ nor $D \hat{\Phi} (X)$ could contract or expand under $\ker \chi$.  As $\chi$ is non-resonant, $D \hat{\Phi} (X)$ would expand or contract under $\ker \chi$ if $D \hat{\Phi} (X)$ was not tangent to the $G$ orbit on $G \times _{\rho(\Gamma)} N/\Lambda$.

In general, $\hat{\phi}$ is of course only continuous.  Instead of analyzing the derivative,  one analyzes the orbit of the unipotent subgroup $\exp t X$ on both suspensions and proves equivariance of $\hat{\Phi}$ under $\exp t X$. 
Then $\hat{\Phi}$ will be $G$ equivariant as the non-resonant roots and the split Cartan $A$ generate $G$ by assumption.

Continuity of $\hat{\Phi}$ now follows from uniqueness arguments.  If $\tau$ is an action by Lipschitz homeomorphisms,  H\"{o}lderness is then standard.



\subsection{Smooth conjugacy}
Let us next outline the ideas behind the proof of Theorem \ref{thm:smoothrigidity}.  Naturally, the authors are searching for a higher rank abelian subgroup $\Sigma$ of $\Gamma$ to which they can apply Theorem \ref{theorem:global rigidity} to.  Then the continuous semi-conjugacy $\Phi$ will automatically be smooth (by uniqueness of semiconjugacies homotopic to the identity).  This will be possible if:
\begin{itemize}
    \item the linear data $\rho (\Sigma)$  contain an Anosov element;
    \item the linear data $\rho (\Sigma)$ do not admit a rank 1 factor.
    \end{itemize}

While existence of an Anosov element in $\Gamma$ is assumed, it may have small centralizer. By works of Prasad and Rapinchuk \cite{Prasad-Rapinchuk2003,MR2150880}, centralizers of generic elements of $\Gamma$ will contain abelian subgroups $\Sigma$  of rank equal to the $\R$-rank of $G$. Such $\Sigma$ will satisfy the requirements above. As the given Anosov element may not be generic, they first show that there is a Zariski dense semigroup  
consisting of Anosov elements.  The proof of this is technical. Essentially they use a  using an invariant cone family, as well as superrigidity (to control the weigths of the linear data) and Ratner's theorem.  

At this point the authors have a regular, $\R$-regular and even hyperregular Anosov element to which they can apply the results of Prasad and Rapinchuk. These also include results on irrationality of the weights of the representation which allows them to exclude rank 1 factors.

\subsection{Comments, problems and related  results} 

Theorem \ref{thm:smoothrigidity} is close to a complete classification of  actions of higher rank lattices on tori and nilmanifolds. One expects this to extend to arbitrary closed manifolds.  However, not much is known.  As mentioned, it would immediately follow if the Anosov Smale conjecture held true.  There are partial results in low dimension compared to the rank of the lattice by Feres and Labourie \cite{Feres95}.  For Anosov actions of semisimple groups, Goetze and Spatzier obtained a smooth classification under a multiplicity free assumption \cite{MR1448015}.  Recently, Damjanovic, Spatzier, Vinhage and Xu \cite{Damjanovic:2022aa} obtained such a  classification  assuming that there is a dense set of Anosov elements in $G$ and that $G$ preserves a volume.  The case of lattices however is much more difficult, and in particular does not follow directly from a suspension argument.
\vspace{.5em}

Finally, we can consider partially hyperbolic actions.  One would again expect considerable rigidity, even classification of such actions.  However, not much is known at this point of time.

\section{Higher rank lattice actions in critical dimensions} \label{sec:latticelowdim}

Zimmer's conjecture concerns actions of higher rank lattices $\Gamma \subset G$ on compact manifolds $M$ of very low dimensions, more precisely in   dimensions so small, smaller than certain \emph{critical  dimensions},  that they preclude the existence of essentially non-trivial actions entirely.   
As mentioned, huge progress on Zimmer's conjecture was made  by Brown, Fisher and Hurtado in 
\cite{Brown-Fisher-Hurtado,MR4132960,Brown:2021wc} even though their results do not yet fully resolve the conjecture.   And for sure, their results    have nothing to say about lattice actions in the critical dimensions and  above.  That is precisely where Zhiren Wang's work \cite{BFHW} with Brown and Rodriguez Hertz comes in.

Let us  be more precise about  the various critical dimensions playing  a role in Zimmer's conjecture.  They only depend on the ambient Lie group $G$ and are defined as follows: 

 \begin{itemize}
    \item $n(G)$  is the smallest dimension for which the ambient group or rather its Lie algebra, has a non-trivial  linear representation.
    \item $v(G)$ is the minimal dimension of a non-trivial homogeneous space $G/H$.  
    \item $d(G)$ is the smallest $v(K)$ where $K$ is a compact form of the complexification of $G$.
\end{itemize}

Then Zimmer's conjecture predicts the following   three phenomena for smooth   $\Gamma$ actions on compact manifolds \cite{Brown:2021wc} (where smoothness can be relaxed to  $C^{1+ \beta}$ for some $\beta$):

\begin{itemize}
    \item  in dimensions below $\min (d(G), v(G))$, $\Gamma$ can only act via finite groups,
    \item in dimensions below $\min (d(G), n(G))$, any volume preserving $\Gamma$ action acts via finite groups,
    \item in dimensions below $n(G)$,  any volume preserving $\Gamma$ action also preserves  a Riemannian metric,
    \item in dimensions below $v(G)$, any $\Gamma$ action has to  preserve   a Riemannian metric.
\end{itemize}

To give a  flavor of their results, \cite{Brown:2021wc} proves the first two cases of Zimmer's conjecture for $\R$-split simple groups.   Actually, they prove a bit more but it becomes  technical rather quickly.

\subsection{Actions near critical dimensions}  Now let us discuss the results in \cite{BFHW}. While Zimmers's conjecture and its known resolution  preclude higher rank lattices from acting in very low dimensions (except via finite groups), we  want to understand and possibly  classify such actions when they do exist.  This is the goal of the Zimmer Program.  It is natural to explore this in the critical range of dimensions just around the conjectural limit for non-existence of actions.    


One sure way of getting faithful actions of a higher rank lattice $\Gamma $ in a semisimple Lie group $G$ is simply by letting $\Gamma$ act on compact homogeneous spaces $G/H$ of $G$ for $H$ a closed subgroup of $G$ (i.e., in dimensions at least $v(G)$). 
When $H$ is connected, then $H$ is a parabolic subgroup of $G$ (essentially by definition).   These give the compact homogeneous spaces of $G$ with the lowest dimension.  Let us call these the {\em boundary actions} or also {\em projective actions} of $\Gamma$ (or $G$).

 As an example, let us briefly describe the parabolic subgroups of $SL(n +1,\R)$.  The smallest (up to conjugacy), called minimal parabolic, consists of all upper triangular matrices (which for $SL(n +1,\R)$ happens to be the maximal solvable algebraic subgroup, also called {\em Borel subgroup}). The other parabolic subgroups consist of matrices which have 0 entries below a suitable ``staircase shape'' (up to conjugacy).  The maximal parabolic subgroups in particular are the stabilizers of the transitive $SL(n,\R)$ action on Grassmannians of $k$-planes, $k=1, \ldots n$.  For them the staircase only has one step, with the zeroes occuring below the $k$'th row and to the left of the $k+1$st column. 

Let us make the crucial remark that boundary actions   behave very differently from other low dimensional actions.  For simplicity, consider $G=SL(n+1,\R)$.  Then $\Gamma =SL(n+1, \Z)$ is a lattice which acts on the $n+1$-torus $\T^{n+1}$ via linear automorphisms.  This is an example of an {\em affine action}.



\medskip

  \textit{\textbf{Fundamental Dichotomy:}} These two types of $\Gamma$ actions are radically different dynamically:
\begin{itemize}
    \item affine actions preserve a probability measure, while
    \item boundary actions never preserve a probability measure.
\end{itemize}

This a key difference.  Indeed, one of the main tools in the Zimmer program, 
 Zimmer's cocycle superrigidity theorem, applies only  when $\Gamma$ preserves a probability measure.  It allows to draw strong measure theoretic conclusions, e.g., positivity of entropy for most elements in $G$ for most actions.  Affine actions fall under this regime, simply because they come from automorphisms which preserve Haar measure.

Boundary actions on the other hand fail to preserve any probability measure.  Indeed,  consider the smallest dimensional boundary action of $PGL(n+1,\R)$   on real projective space $\RP ^{n}$.  There  always exist elements $\gamma \in \Gamma$ which are diagonalizable with distinct eigenvalues (by work of  Prasad and  
Raghunathan \cite{Prasad-Raghunathan1972}).   Such $\gamma$  have finitely many fixed points acting on 
$\RP ^n$ to which all other points get attracted (they are Morse-Smale diffeomorphisms, generalizing the well-known North-South dynamics from boundary actions of rank 1 lattices).  Hence any $\gamma$ invariant probability measure is atomic.  Therefore  $\Gamma$-invariant probability measures are also atomic and thus $\Gamma$ will have a periodic orbit which contradicts Zariski density of $\Gamma$ in $G$.  This argument generalizes to all boundary actions of $\Gamma$ on any boundary $G/Q$ of $G$, for $Q$ a parabolic subgroup of $G$.

 Next consider  actions of lattices in $SL(n,\R)$ for simplicity. The dichotomy results in the main theorem of \cite{BFHW}:

\begin{theorem}
    Suppose $\Gamma$ is a uniform lattice in $SL(n+1, \R)$ with a smooth action $\alpha$  on a  compact manifold $M$ with $\dim M <n$.  Then $\Gamma$  fixes a Borel probability measure on $M$.  If $\dim M =n$, then either there exists a $\Gamma$-invariant probability measure or the action is measurable conjugate to a boundary action of  $\Gamma$  on $G/Q$ for some parabolic $Q \subset G$. 
\end{theorem}

The authors recently sharpened this theorem  to get smooth conjugacies, cf. Section \ref{sec5:comments}. 
We also note that the first claim in this theorem was first proved in \cite{Brown:2021wc}.  However, Wang and co-authors prove a more general new result in \cite{BFHW} as we shall see. 


Again, as a main tool for their results, the authors  analyze the suspension action from $\Gamma$ to $G= SL(n+1, \R)$  (cf. \ref{sec:susp}).  Due to problems with escape of mass, the case of non-uniform lattices again gets much more technical, and we will only consider actions of uniform lattices here. It is rather impressive that they can manage these difficulties in the non-uniform case.

Resonance will again play a major role in their considerations, as in Section \ref{sec:semiconj-nr}.  We will need the following definitions to state their results in full, for general semisimple $G$ as above.

 Let $\mathfrak g$ be a semisimple Lie algebra over $\R$, and let $\Sigma$ be it set of restricted roots relative to a maximally split Cartan subalgebra. Let $\Sigma ^+$ denote the set of positive roots (for some choice of ordering for the roots or equivalently a choice of positive Weyl chamber in the split Cartan).  Suppose $\mathfrak q$ is a parabolic subalgebra of $\mathfrak g$, i.e., the Lie algebra of a parabolic subgroup. By their structure theory,  every parabolic subalgebra is conjugate to one of finitely many standard ones which are sums of certain root spaces for the split Cartan picked before. Note the following:   if a restricted root $\alpha$ has its root space $\mathfrak g ^{\alpha} $ in $\mathfrak q$,  so does every positive multiple of $\alpha$ in $\Sigma$. 
 
Denote by $\hat{\Sigma}$ the set of   (restricted) roots up to positive scalar multiple.  For $\sigma \in \hat{\Sigma}$, we call the sum of the root spaces $\mathfrak g _{\sigma} = \oplus _{\alpha \in \sigma} \mathfrak g ^{\alpha}$ a {\em coarse (restricted) root space}. 
Then each $\mathfrak g _{\sigma}$
is contained in a Borel subgroup.
 Indeed,  if a restricted root $\alpha$ has its root space $\mathfrak g ^{\alpha} $ in some parabolic subalgebra $\mathfrak q$,  so does every positive multiple of $\alpha$. This comment allows us to make the following definition.
   
\begin{Definition}
    The \emph{resonant codimension} of $\mathfrak q$ is the cardinality of 
   $ \{ \alpha \in \hat{\Sigma}\mid \mathfrak g ^{\alpha } \not\subset \mathfrak{q} \}$.
 We let the {\em minimal resonant codimension} $r(\mathfrak g)$ of $\mathfrak g$ be the smallest resonant codimension of any proper parabolic subalgebra. If $G$ is a semisimple Lie group with Lie algebra $\mathfrak g$, we set $r(G) = r(\mathfrak g)$.
\end{Definition}

For $SL(n+1,\R)$, the resonant dimension is $n$ as is easily seen from checking all maximal parabolics, i.e., the ones with Grassmannians as quotient spaces. The authors explicitly determine the resonant dimensions of many classical groups. For real split simple groups $G$, it is easily seen that $r(G)=v(G)$ as all root spaces are one-dimensional and distinct roots are never positively proportional.

Now we can state their main technical result for actions on manifolds of near critical dimensions:

\begin{theorem}  \label{them:low}
    Suppose a higher rank lattice $\Gamma \subset G$ acts smoothly ($C^{1 +\theta}$) on  a compact manifold $M$.  
    \begin{itemize}
   \item  If  $\dim M < r(G)$, then $\Gamma$ preserves a Borel probability measure. 
   \item If $\dim M = r(G)$, then 
   \begin{itemize}
   \item either $\Gamma$ preserves a Borel probability measure, or
   \item  there is a parabolic $Q \subset G$ and a $\Gamma$ quasi-invariant probability measure on $M$ such that the $\Gamma$ action on $M$  is measurably conjugate to a finite extension of the boundary action $\Gamma$ on  $G/Q$. 
    \end{itemize}
     \end{itemize}
\end{theorem}

\subsection{Intermediate dimensions}

When the dimension of $M$ is higher than the critical dimensions,  one can construct many actions by combining projective,  affine and trivial actions, and also suspend actions of parabolic subgroups of $G$.  Examples of the latter include {\em homogeneous} actions given by left translations of $\Gamma$ on a compact homogeneous space  via an embedding $\rho$  of $\Gamma$ into $H$. Furthermore, if $\rho (\Gamma) $ normalizes a closed subgroup $M \subset H$, then we get a  $\Gamma$ action on $X= M \setminus H/\Lambda$.  We call the latter actions {\em bi-homogenous} 
if $X$ is a manifold. Notably, there are  bi-homogeneous examples that are not  homogeneous (even for $M$ compact), cf.     \cite[Example 4.9]{Damjanovic:2023aa}  
(that this example is not homogeneous follows from considerations about the fundamental group). When $M$ is not compact, there are further constructions of homogeneous quotients which cover compact manifolds \cite{Margulis-problems, MR816674,MR1410076}.  
More subtly, one can construct actions by ``blowing up'' an affine action over  fixed points or periodic orbits,  as was first done by Katok and Lewis, then Benveniste \cite{KL,Benveniste-Fisher}. Recently, Fisher and Melnick constructed new examples even in low dimensions \cite{Fisher:2022aa}.

One can vaguely conjecture that  all higher rank  smooth $\Gamma$ actions arise from such constructions. Making this more  precise is not very clear at this point.  If the action preserves a volume and is either ergodic or transitive, then 
conjecturally the action fibers over   a $\Gamma$ action with bi-homogeneous orbit closures such that the factor map restricts to a diffeomorphism on a  dense open set. 
This is similar to conjectures of Zimmer \cite{Zimmer-ICM} and Margulis \cite[Problem 11]{Margulis-problems}.  
However, these conjectures are not approachable using current methods.   

Even more generally, if an action $\alpha$ acts only absolutely continuously w.r.t. a volume, 
as in the case of projective actions, 
one can hope that $\alpha$  factors over a projective action and ``preserves'' (conditional) measures along the fibers.  
In the measurable category, the work of Nevo and Zimmer achieves such a structure theorem under additional  assumptions 
\cite{NevoZ00,NevoZ02}. Unfortunately, those are difficult to verify. 

Enter once again Brown, Rodriguez Hertz and Wang who bring these hopes  closer to reality,  proving a  very satisfying  structure theorem for actions in ``sufficiently  low'' dimensions.  First, to make this precise, they define another invariant, $m(G)$ of $G$, defined in order to  capture the dimensions between quotients of $G$ by maximal parabolics and by next to maximal ones.   The goal is to force an action in this range of dimensions to either preserve a measure or to factor over  $G/Q$ where $Q$ is a maximal parabolic. 

\begin{Definition}
  Let  $m (\mathfrak g)$ and $m(G)$ denote the minimal resonant codimension of non-maximal proper parabolic subalgebras of $\mathfrak g$. 
\end{Definition}

 As an example, one can easily calculate that $m(SL(n+1,\R)) = 2n-1$ by considering partial flag varieties.

Second, they carefully define a suitable notion of factorization:
Suppose we are given an action
$\alpha$  on a measure space $(X,\nu)$.  Suppose  $(Y, \eta)$ is another measure space.  Let $\operatorname{Aut}(Y, \eta)$ denote the group of measure preserving invertible transformation of $Y$.  

Suppose  $\psi: \Gamma \times X \mapsto \operatorname{Aut}(Y,\eta)$ is a cocycle over $\alpha$ (i.e., satisfies the ``chain rule'').  Then we can build the extension $\alpha _{\psi}$ of $\alpha $ by $\psi$ on $X \times Y$ by setting
    \[\gamma (x,y) = (\alpha (\gamma ) (x), \psi (\gamma, x) (y)).\]

\begin{Definition}
    Call an  action  $\beta$  of $\Gamma $ on a measure space $(Z, \mu)$    a {\em relatively measure-preserving extension} if $\beta$ is measure theoretically equivalent to an extension $\alpha _{\psi}$ for suitable $\alpha$ and $\psi$.
    \end{Definition}

    Finally, we can state their last, marvelous result:

    \begin{theorem}
        If a higher rank lattice $\Gamma \subset G$ acts on a compact manifold $M$ with  $\dim M < m(G)$, then, for some $\Gamma$ quasi-invariant measure on $M$,  the action is a relatively measure-preserving extension of a projective action.
    \end{theorem}

In particular, for lattices in $SL(n+1,\R)$, with  $n\geq 2$, all smooth actions on compact manifolds of dimension less than $2n-1$ are measure preserving extensions of boundary actions. Beware that the latter could be the trivial action, i.e., $\beta$ preserves a probability measure.    Thus, the lack of an   invariant measure
is fully due to the presence of a projective factor.

\subsection{Outline of the arguments}

Given an action $\alpha$ of a higher rank lattice $\Gamma$ on a compact manifold $M$, the authors first pass to the suspension $\tilde{\alpha}$ on the suspension space  $M^{\alpha}$ - again because this allows them to argue with actions of connected Lie groups rather than discrete groups.   Resonance will play an essential role again, especially in forcing invariance properties of measures.  This time, however, resonance is not with respect to linear data as in Section \ref{sec:anosovlattice}.    Indeed, $M$ a priori is an arbitrary compact manifold, and linear data are simply not defined.  Instead, they use Lyapunov exponents $\lambda$ of the action of the maximal split Cartan  $A \subset G$ acting on $M^{\alpha}$, with respect to an $A$-invariant probability measure.
 These are linear functionals on $A$ (identified with $\R ^k$) whose values on an element $a \in A$  are simply the Lyapunov exponents of the diffeomorphism on $M^{\alpha}$ by determined $a$.    We will call these linear functionals  {\em Lyapunov exponents}. 
 This is in the same spirit as how roots and weights are defined for linear representations of $G$. These definitions can be made for any action of  $\R^k$ on a compact manifold w.r.t an invariant probability measure.
 
 In the presence of an invariant probability measure, one can define stable distributions for particular elements of $A$ and refine the usual Oseledets decomposition into a joint decomposition for all elements of $A$. 
 More useful is the related decomposition into {\em coarse Lyapunov spaces} where one bundles positively proportional Lyapunov functionals together, and thinks of them as {\em coarse Lyapunov exponents}.  The sum of all the Oseledets spaces for all coarsely related Lyapunov exponents is called a {\em coarse Lyapunov space}.  They determine the {\em  coarse Lyapunov distributions} which are integrable into measurable sub-foliations of the various stable foliations.  We call them {\em coarse Lyapunov foliations}. 
 This has been all carefully proved and explained in \cite{Brown-Hertz-Wang2023-JMD} and is summarized in Section 3 of their paper.  While they develop the needed notions and estimates for both uniform and non-uniform lattices, we will only discuss the case of uniform lattices $\Gamma$,  for simplicity.

\vspace*{1em}

\emph{Measure rigidity ideas and non-resonance.}  The main goal is to show invariance of certain measures $\mu$ on $M^{\alpha}$ under the suspension action of  $G$  if they are invariant under a split Cartan $A \subset G$.  This is similar to the problem of classifying probability measures on homogeneous spaces $G/\Gamma$ which are invariant under a split Cartan subgroup $A$, the so-called {\em measure rigidity problem}.  Again, for higher rank $G$,  $A$ invariance should imply invariance under a much larger subgroup or even the full group  $G$.  This has been much investigated over the last three decades.  Indeed,   Einsiedler and Lindenstrauss' work \cite{EinLin2015}  is the culmination of  current methods.


In their paper, Brown, Rodriguez Hertz and Wang adopt these measure rigidity ideas to their end.  First note that it suffices to prove invariance of $\mu$ under the unipotent root subgroups $U$ as these generate $G$ as a Lie group.




By a long, complicated argument they show that the conditional measures of $\mu$ of on the orbit foliations of the unipotent root subgroups $U \subset G$ are invariant under $U$.    Since these $U$ generate $G$, this will suffice.

More precisely, consider probability measures $\mu$ on $M^{\alpha}$ which project to the Haar measure $\nu$ on $G/\Gamma$ and are invariant under a minimal parabolic subgroup $P \subset G$.  Since $P$ is amenable such measures always exist by averaging probability measures on $M^{\alpha}$ over F\o lner sets in $P$.  Consider  $Q$  the stabilizer of $\mu$ in $G$.   Since $P \subset Q$, $Q$ is a parabolic subgroup.  The goal is to show that $Q=G$. 

Henceforth, given a root $\beta$ of $G$,  let $U$ be the unipotent subgroup generated by the unipotent subgroups corresponding to all positively proportional roots $c \beta, c >0$.  We will also call $U$ the coarse root subgroup determined by $\beta$.  Since such $U$ generate all of $G$, it suffices to show invariance of $\mu$ under such $U$.  For this, it suffices that the conditional measures of $\mu$ for the $U$ orbit foliations are absolutely continuous along $U$ orbits.  

Note that the $U$ orbit foliation is a smooth subfoliation of the coarse Lyapunov foliation for the coarse Lyapunov exponent determined by $\beta$. Hence $\mu$ will have absolute continuous conditional measures along the $U$-orbits if $\mu$ has absolutely continuous conditional measures for the coarse foliation determined by $\beta$.

 Here entropy enters.  Indeed, Pesin famously proved that the measure-theo-retic entropy of a smooth diffeomorphism with respect to an invariant absolutely continuous probability measure  as a sum of the positive Lyapunov exponents (assuming ergodicity). Ledrappier and Young developed a general entropy formula in  \cite{ledrappier_young1985}  which gave a converse to Pesin's work and forced absolute continuity of conditional measures on stable manifolds in the equality case.  Before, Ledrappier had a similar result for conditional entropy along measurable foliations with smooth leaves \cite{Led84}.
 
In consequence, they need to show that the conditional entropy of $\mu$ along the coarse Lyapunov foliations of roots is given by the sum of the exponents times the dimension of the leaves of a coarse Lyapunov foliation on $M^{\alpha}$.    We are looking for a simple criterion on the root $\beta$  which guarantees this. Note that  $U$ is expanded and contracted by $A$ according to the root $\beta$ defining $U$. Let $\hat{\beta}$  be  the coarse Lyapunov functional determined by $\beta$ on $G/\Gamma$ and  $\chi$ be the one on $M^{\alpha}$, again  determined by $\beta$. Let $W^{\chi}$ be its coarse Lyapunov foliation.  Since $M^{\alpha}$ fibers $G$-equivariantly over $G/\Gamma$, an Abramov--Rohlin type  formula relates conditional entropy for $W^{\chi}$ with the base entropy for the $W^{\beta}$ foliation on $G/\Gamma$ and the fiberwise entropy of $W^{\chi}$.  However, if $\beta$ is {\bf non-resonant } with any fiber weight, then there is no contribution from the fiber, and the entropies for the $W^{\chi}$ foliation on $M^{\alpha}$ and the $W^{\beta}$ foliation on $G/\Gamma$ coincide.  From this one gets invariance of $\mu$ under $U$.


\smallskip

{\em How to apply non-resonance?}
Consider first the simple case when the dimension of the manifold is less than the resonant codimension of $G$:  $\dim M < r(G)$.  Then there are at most $r(G) -1$ resonant Lyapunov exponents of $\mu$ coming from the fiber of the projection $M^{\alpha}\rightarrow G/\Gamma$.  Thus the resonant co-dimension of the stabilizer $Q$ of $\mu$ is at most $r(G)-1$.  Hence $Q=G$ by definition of the resonant co-dimension.

Now suppose the dimension of $M$ is $r(G)$.  If there is a non-resonant fiberwise Lyapunov exponent then there are at most $r(G)$-1 coarse resonant roots, and we can run the same argument as when  $\dim M < r(G)$.  Thus we conclude that the $P$ -invariant measure $\mu$ is $G$-invariant.  

In consequence,  we may assume that all fiberwise Lyapunov exponents of $\mu$ are resonant with a root.  Since $P$ leaves $\mu$ invariant, we may assume that the latter roots are negative roots.  
Then some $s \in A$ will contract all fibers, and it follows quickly that the fiber measures of $\mu$ have finitely many atoms which will be 
fixed by a parabolic subgroup $Q$. They end up viewing the fiber as a finite extension of $G/Q$ using the finiteness of the atoms. 

The case of intermediate dimensions is considerably more difficult and technical, and we will refrain from discussing details. As a key additional argument, they are using the so-called 
Einsiedler--Katok lemma from measure rigidity  which guarantees strong invariance properties of the conditional measures $\mu _x$ of $\mu$ along  suitable homogeneous directions in the coarse Lyapunov foliations. More precisely, the $\mu _x$ will be supported on orbits of subgroups $S_x$ and bi-invariant under suitable subgroups $H_x \subset S_x$.

\vspace*{-0.02em}
\subsection{Comments, problems, and related results}\label{sec5:comments}

\vspace*{-0.41em}
\subsubsection{Factor theorems}
The factorization  results from \cite{BFHW} discussed above analyze smooth actions of higher rank lattices $\Gamma$, possibly factoring to a projective action of $\Gamma$, or boundary action on $G/Q$, $Q$ a parabolic more generally.  This stands in direct contrast with the celebrated work of Margulis on measurable quotients of boundary  actions of $\Gamma$.  Margulis shows that they are all  measurably isomorphic to some (potentially different) boundary action \cite{Margulis-book1991, Zimmer1}. No smoothness of the factor action is assumed for the latter.  Margulis used it to prove his famous finiteness theorem that normal subgroups of irreducible higher-rank lattices are either central or of finite index.  While pulling in opposite directions, there are commonalities in these results which allowed Brown, Rodriguez Hertz and Wang to give a different proof of the measurable factor theorem in \cite{MR4422052}. They apply ideas from higher rank measure rigidity which they also so centrally use in \cite{BFHW}.  

For continuous factors of projective actions, Dani proved in a similar vein that such factors are $C^0$ conjugate to another projective action \cite{Dani}.  A related smooth version was proved  in \cite{GS-smoothfactors}, resulting in a  smooth conjugacy.  However, they needed existence of ``differentiable sinks'', i.e., fixed points of at least one suitable element $\gamma \in \Gamma$ whose derivatives have all eigenvalues less than 1 in norm. In low dimensions, e.g., $n-1$ for lattices in $SL(n,\R)$,  such results immediately follow  from the forthcoming smooth versions of factorization results in \cite{BFHW}, at least for real split semisimple groups $G$ (without assuming existence of differentiable sinks).

\vspace*{-0.2em}
\subsubsection{Local rigidity}
Given an action of a finitely generated group, e.g., any higher rank lattice, on a compact manifold $M$, we call it {\em locally rigid} if a perturbation of the action  close on a finite generating set is $C^{k}$-conjugate to the original action, for some $k$. We can measure the closeness in terms of a $C^l$ topology, for some $l$.  Thus we can speak of $C^{k,l}$-local rigidity.  In addition, we can require the conjugacy to be close to the identity.  

Local rigidity was proved in the 1990s for projective actions of higher rank lattices: for the $C^{4,\infty}$ case by Kanai \cite{Kanai}  for uniform lattices in $SL(n,\R)$, using geometric means, and by Katok and Spatzier in \cite{KS97} for any uniform higher rank lattice for the $C^{1,\infty}$ case.  The latter used a mix of dynamical and representation theoretic methods.   Naturally, the classification results of \cite{BFHW} automatically give local rigidity in low dimensions, e.g., for  lattices in $SL(n+1,\R)$ on $\RP ^n$, for $n>1$. Notably, these are the first results when  the lattice is not co-compact. 

Remarkably, $C^{0,0}$-local rigidity has been studied for many actions, starting with Sullivan for actions of convex cocompact Kleinian groups on their limit sets in \cite{MR806415}.   One twist here is that one can only hope to get a semi-conjugacy from the perturbation to the original action, i.e., we get a surjective equivariant continuous map which may however not be invertible.  For actions of Gromov hyperbolic groups on their boundaries, this was recently proved by Mann, Manning, and Weisman \cite{Mann:2022aa}, following earlier works by Bowden and Mann \cite{MR4468857} as well as Mann and Manning \cite{Mann:2021aa}. 

For projective actions of uniform higher rank lattices, $C^{0,0}$-local rigidity was achieved by Connell, Islam, Nguyen, and  Spatzier in \cite{Connell:2023aa}, using quasi-isometric rigidity and Morse Lemma type techniques. Earlier results for higher-rank lattices were obtained in the Lipschitz category by Kapovich, Kim, and Lee \cite{kapovich_kim_lee_2022}, using Sullivan type arguments.

\subsection*{Acknowledgments}  We thank David Fisher for several discussions about the Zimmer program. We also thank the referee whose careful reading of the text greatly improved its readability.

\bibliographystyle{numeric}

\begin{thebibliography}{99}
\bibitem{Barthelme:2021aa} 
\newblock T. Barthelm\'{e}, C. Bonatti, A. Gogolev and F. Rodriguez~Hertz,
\newblock {Anomalous {A}nosov flows revisited},
\newblock {\em Proc. Lond. Math. Soc. (3)}, \textbf{122} (2021), 93--117.

\bibitem{Benveniste-Fisher} 
\newblock E.~J. Benveniste and D. Fisher,
\newblock {Nonexistence of invariant rigid structures and invariant almost rigid structures},
\newblock {\em Comm. Anal. Geom.}, \textbf{13} (2005), 89--111.

\bibitem{Bonatti:2023aa} 
\newblock C. Bonatti and I. Iakovoglou,
\newblock {Anosov flows on 3-manifolds: The surgeries and the foliations},
\newblock {\em Ergodic Theory Dynam. Systems}, \textbf{43} (2023), 1129--1188.

\bibitem{MR4468857} 
\newblock J. Bowden and K. Mann,
\newblock {{$C^0$} stability of boundary actions and inequivalent {A}nosov flows},
\newblock {\em Ann. Sci. \'{E}c. Norm. Sup\'{e}r. (4)}, \textbf{55} (2022), 1003--1046.

\bibitem{Brown-Fisher-Hurtado} 
\newblock A. Brown, D. Fisher and S. Hurtado,
\newblock Zimmer's conjecture: Subexponential growth, measure rigidity, and strong property (T),
\newblock \emph{Ann. of Math. (2)}, \textbf{196} (2022),  891--940.

\bibitem{MR4132960} 
\newblock A. Brown, D. Fisher and S. Hurtado,
\newblock {Zimmer's conjecture for actions of {${\rm SL}(m, \Bbb Z)$}},
\newblock {\em Invent. Math.}, \textbf{221} (2020), 1001--1060.

\bibitem{Brown:2021wc}
\newblock A. Brown, D. Fisher and S. Hurtado,
\newblock Zimmer's conjecture for non-uniform lattices and escape of mass,
\newblock preprint, \arxiv{2105.14541}, 2021.


\bibitem{MR3702679} 
\newblock A. Brown, F. Rodriguez~Hertz and Z. Wang,
\newblock {Global smooth and topological rigidity of hyperbolic lattice actions},
\newblock {\em Ann. of Math. (2)}, \textbf{186} (2017), 913--972.

\bibitem{BFHW} 
\newblock A. Brown, F. Rodriguez~Hertz and Z. Wang,
\newblock {Invariant measures and measurable projective factors for actions of higher-rank lattices on manifolds},
\newblock {\em Ann. of Math. (2)}, \textbf{196} (2022), 941--981.

\bibitem{MR4422052}
\newblock A. Brown, F. Rodriguez~Hertz and Z. Wang,
\newblock {The normal subgroup theorem through measure rigidity},
\newblock in {\em Dynamics, Geometry, Number Theory---the Impact of {M}argulis on Modern Mathematics}, 66--91, Univ. Chicago Press, Chicago, IL, 2022.

\bibitem{Brown-Hertz-Wang2023-JMD}
\newblock A. Brown, F. Rodriguez~Hertz and Z. Wang,
\newblock {Smooth ergodic theory of $\Bbb Z^d$-actions}, 
\newblock {\em J. Mod. Dyn.}, {\bf 19} (2023), 455--540.

\bibitem{MR1911660}
\newblock M.~Burger and N.~Monod,
\newblock {Continuous bounded cohomology and applications to rigidity theory},
\newblock {\em Geom. Funct. Anal.}, \textbf{12} (2002), 219--280.

\bibitem{Connell:2023aa}
\newblock C. Connell, M. Islam, T. Nguyen and R. Spatzier,
\newblock Boundary actions of lattices and $c^0$ local semi-rigidity,
\newblock preprint, \arxiv{2303.00543}, 2023.

\bibitem{Damjanovic:2022aa}
\newblock D. Damjanovic, R. Spatzier, K. Vinhage and D. Xu,
\newblock Anosov actions: classification and the {Z}immer program,
\newblock preprint, \arxiv{2211.0819511}, 2022.

\bibitem{Damjanovic:2019aa}
\newblock D. Damjanovic, A. Wilkinson and D. Xu,
\newblock {Pathology and asymmetry: Centralizer rigidity for partially hyperbolic diffeomorphisms},
\newblock \emph{Duke Math. J.}, \textbf{170} (2021),  3815--3890.

\bibitem{Damjanovic:2023aa} 
\newblock D. Damjanovic, A. Wilkinson and D. Xu,
\newblock {Transitive centralizers and fibered partially hyperbolic systems},
\newblock \emph{Int. Math. Res. Not. IMRN}, (2024),  9686--9704.

\bibitem{MR4145817}
\newblock D. Damjanovi\'{c} and D. Xu,
\newblock {On classification of higher rank {A}nosov actions on compact manifold},
\newblock {\em Israel J. Math.}, \textbf{238} (2020), 745--806.

\bibitem{MR518242} 
\newblock S.~G. Dani,
\newblock {Nilmanifolds with {A}nosov automorphism},
\newblock {\em J. London Math. Soc. (2)}, \textbf{18} (1978), 553--559.

\bibitem{Dani} 
\newblock S.~G. Dani,
\newblock {Continuous equivariant images of lattice-actions on boundaries},
\newblock {\em Ann. of Math. (2)}, \textbf{119} (1984), 111--119.

\bibitem{MR2140439}
\newblock S.~G. Dani and Meera~G. Mainkar,
\newblock {Anosov automorphisms on compact nilmanifolds associated with graphs},
\newblock {\em Trans. Amer. Math. Soc.}, \textbf{357} (2005), 2235--2251.

\bibitem{Dere:2020aa}
\newblock J. Der{\'e} and M. Mainkar,
\newblock {Anosov diffeomorphisms on infra-nilmanifolds associated to graphs},
\newblock {\em Math. Nachr.},  \textbf{296} (2023), 610--629.

\bibitem{EinLin2015}
\newblock M. Einsiedler and E. Lindenstrauss,
\newblock {On measures invariant under tori on quotients of semisimple groups},
\newblock {\em Ann. of Math. (2)}, \textbf{181} (2015), 993--1031.

\bibitem{Feres95}
\newblock R.~Feres and F.~Labourie,
\newblock {Topological superrigidity and {A}nosov actions of lattices},
\newblock {\em Ann. Sci. {\'E}cole Norm. Sup. (4)}, \textbf{31} (1998), 599--629.

\bibitem{MR4521507}
\newblock D. Fisher,
\newblock Groups acting on manifolds: Around the {Z}immer program,
\newblock in {\em Group Actions in Ergodic Theory, Geometry, and Topology--Selected Papers},  609--683, Univ. Chicago Press, Chicago, IL, 2020.

\bibitem{FKS}
\newblock D. Fisher, B. Kalinin and R. Spatzier,
\newblock {Global rigidity of higher rank {A}nosov actions on tori and nilmanifolds}, with an appendix by James F. Davis,
\newblock {\em J. Amer. Math. Soc.}, \textbf{26} (2013), 167--198.

\bibitem{Fisher:2022aa} 
\newblock D. Fisher and K. Melnick,
\newblock {Smooth and analytic actions of $SL(n,{\bf R})$ and $SL(n,{\bf Z})$ on closed $n$-dimensional manifolds},
\newblock \emph{Kyoto J. Math.}, \textbf{64} (2024), 873--904.

\bibitem{Franks1970}
\newblock J. Franks,
\newblock {Anosov diffeomorphisms},
\newblock in {\em Global {A}nalysis ({P}roc. {S}ympos. {P}ure {M}ath., {V}ol. {{XIV}}, {B}erkeley, {C}alif., 1968)}, 61--93, Amer. Math. Soc., Providence, RI, 1970. 

\bibitem{MR2219247} 
\newblock J. Franks and M. Handel,
\newblock {Distortion elements in group actions on surfaces},
\newblock {\em Duke Math. J.}, \textbf{131} (2006), 441--468.

\bibitem{MR1259430}
\newblock \'{E}. Ghys,
\newblock {Rigidit\'{e} diff\'{e}rentiable des groupes fuchsiens},
\newblock {\em Inst. Hautes \'{E}tudes Sci. Publ. Math.}, \textbf{78} (1993), 163--185.

\bibitem{MR1448015}
\newblock E.~R. Goetze and R.~J. Spatzier,
\newblock {On {L}iv\v{s}ic's theorem, superrigidity, and {A}nosov actions of semisimple {L}ie groups},
\newblock {\em Duke Math. J.}, \textbf{88} (1997), 1--27.

\bibitem{MR816674}
\newblock W.~M. Goldman,
\newblock {Nonstandard {L}orentz space forms},
\newblock {\em J. Differential Geom.}, \textbf{21} (1985), 301--308.

\bibitem{MR3413978}
\newblock A. Gorodnik and R. Spatzier,
\newblock {Mixing properties of commuting nilmanifold automorphisms},
\newblock {\em Acta Math.}, \textbf{215} (2015), 127--159.

\bibitem{GS-smoothfactors}
\newblock A. Gorodnik and R. Spatzier,
\newblock {Smooth factors of projective actions of higher-rank lattices and rigidity},
\newblock {\em Geom. Topol.}, \textbf{22} (2018), 1227--1266.

\bibitem{GreenTao}
\newblock B. Green and T. Tao,
\newblock {The quantitative behaviour of polynomial orbits on nilmanifolds},
\newblock {\em Ann. of Math. (2)}, \textbf{175} (2012), 465--540.

\bibitem{GromovPolyGrowth}
\newblock M. Gromov,
\newblock {Groups of polynomial growth and expanding maps},
\newblock {\em Inst. Hautes {\'E}tudes Sci. Publ. Math.}, \textbf{53} (1981), 53--73.

\bibitem{Hurder}
\newblock S. Hurder,
\newblock {Rigidity for {A}nosov actions of higher rank lattices},
\newblock {\em Ann. of Math. (2)}, \textbf{135} (1992), 361--410.

\bibitem{hurder-survey1994}
\newblock S. Hurder,
\newblock {A survey of rigidity theory for {A}nosov actions},
\newblock in {\em Differential Topology, Foliations, and Group Actions ({R}io de {J}aneiro, 1992)},  143--173, {Contemp. Math.}, vol. 161, Amer. Math. Soc., Providence, RI, 1994.

\bibitem{MR2240907}
\newblock B. Kalinin and V. Sadovskaya,
\newblock {Global rigidity for totally nonsymplectic {A}nosov {$\Bbb Z^k$} actions},
\newblock {\em Geom. Topol.}, \textbf{10} (2006), 929--954.

\bibitem{KalSad07}
\newblock B. Kalinin and V. Sadovskaya,
\newblock {On the classification of resonance-free {A}nosov {$Z^k$} actions},
\newblock {\em Michigan Math. J.}, \textbf{55} (2007), 651--670.

\bibitem{KaSp04}
\newblock B. Kalinin and R. Spatzier,
\newblock {On the classification of {C}artan actions},
\newblock {\em Geom. Funct. Anal.}, \textbf{17} (2007), 468--490.

\bibitem{Kanai}
\newblock M.~Kanai,
\newblock {A new approach to the rigidity of discrete group actions},
\newblock {\em Geom. Funct. Anal.}, \textbf{6} (1996), 943--1056.

 \bibitem{kapovich_kim_lee_2022} 
\newblock M. Kapovich, S. Kim and J. Lee,
\newblock {Structural stability of meandering-hyperbolic group actions},
\newblock {\em J. Inst. Math. Jussieu}, \textbf{23} (2024),   753--810.

\bibitem{KatokHasselBook}
\newblock A. Katok and B. Hasselblatt,
\newblock {\em Introduction to the Modern Theory of Dynamical Systems},
\newblock Cambridge University Press, Cambridge, 1995. 

\bibitem{KL1}
\newblock A.~Katok and J.~Lewis,
\newblock {Local rigidity for certain groups of toral automorphisms},
\newblock {\em Israel J. Math.}, \textbf{75} (1991), 203--241.

\bibitem{KL}
\newblock A.~Katok and J.~Lewis,
\newblock {Global rigidity results for lattice actions on tori and new examples of volume-preserving actions},
\newblock {\em Israel J. Math.}, \textbf{93} (1996), 253--280.

\bibitem{KL2}
\newblock A.~Katok, J.~Lewis and R.~Zimmer,
\newblock {Cocycle superrigidity and rigidity for lattice actions on tori},
\newblock {\em Topology}, \textbf{35} (1996), 27--38.

\bibitem{KS97}
\newblock A.~Katok and R.~J. Spatzier,
\newblock Differential rigidity of {A}nosov actions of higher rank abelian groups and algebraic lattice actions,
\newblock {\em Proc. Steklov Inst. Math.}, \textbf{216} (1997), 287--314.

\bibitem{MR1410076}
\newblock F.~Labourie,
\newblock Quelques r\'{e}sultats r\'{e}cents sur les espaces localement homog{\`e}nes compacts,
\newblock in {\em Manifolds and Geometry ({P}isa, 1993)},  267--283, Sympos. Math., vol. XXXVI, Cambridge Univ. Press, Cambridge, 1996.

\bibitem{MR2471923}
\newblock J. Lauret and C.~E. Will,
\newblock {Nilmanifolds of dimension {$\leq 8$} admitting {A}nosov diffeomorphisms},
\newblock {\em Trans. Amer. Math. Soc.}, \textbf{361} (2009), 2377--2395.

\bibitem {Led84}
\newblock F. Ledrappier,
\newblock {Propri\'et\'es ergodiques des mesures de Sina\"i}, 
\newblock {\em Inst. Hautes \'Etudes Sci. Publ. Math.},  \textbf{59} (1984), 163--188.


\bibitem{ledrappier_young1985}
\newblock F. Ledrappier and L.-S. Young,
\newblock {The metric entropy of diffeomorphisms: Part II: Relations between entropy, exponents and dimension},
\newblock {\em Ann. of Math. (2)}, \textbf{122} (1985), 540--574.

\bibitem{Malcev:1951aa}
\newblock A.~I. Malcev,
\newblock On a class of homogeneous spaces,
\newblock {\em Amer. Math. Soc. Translation}, (1951), no.~39, 33 pp.

\bibitem{Mane77}
\newblock R. Ma{\~n}{{\'e}},
\newblock {Quasi-{A}nosov diffeomorphisms and hyperbolic manifolds},
\newblock {\em Trans. Amer. Math. Soc.}, \textbf{229} (1977), 351--370.

\bibitem{Mann:2021aa}
\newblock K. Mann and J.~F. Manning,
\newblock {Stability for hyperbolic groups acting on boundary spheres},
\newblock \emph{Forum Math. Sigma}, \textbf{11} (2023), Paper No. e83, 25 pp.

\bibitem{Mann:2022aa}
\newblock K. Mann, J.~F. Manning and T. Weisman,
\newblock Stability of hyperbolic groups acting on their boundaries,
\newblock preprint, \arxiv{2206.14914}, 2022.

\bibitem{Manning74}
\newblock A. Manning,
\newblock {There are no new {A}nosov diffeomorphisms on tori},
\newblock {\em Amer. J. Math.}, \textbf{96} (1974), 422--429.

\bibitem{Margulis-problems}
\newblock G. Margulis,
\newblock Problems and conjectures in rigidity theory,
\newblock in {\em Mathematics: Frontiers and Perspectives}, 161--174, Amer. Math. Soc., Providence, RI, 2000. 

\bibitem{Margulis-ICM1975}
\newblock G.~A. Margulis,
\newblock Discrete groups of motions of manifolds of nonpositive curvature,
\newblock in {\em Proceedings of the {I}nternational {C}ongress of {M}athematicians ({V}ancouver, {B}.{C}., 1974), Vol. 2}, 21--34, Canad. Math. Congress, Montreal, QC, 1975.

\bibitem{margulis-book}
\newblock G.~A. Margulis,
\newblock {Finiteness of quotient groups of discrete subgroups},
\newblock {\em Funktsional. Anal. i Prilozhen.}, \textbf{13} (1979), 28--39.

\bibitem{Margulis-book1991}
\newblock G.~A. Margulis,
\newblock {\em Discrete Subgroups of Semisimple {L}ie Groups}, {Ergeb. Math. Grenzgeb. (3)}, vol. 17, 
\newblock Springer-Verlag, Berlin, 1991.

\bibitem{NevoZ00}
\newblock A. Nevo and R.~J. Zimmer,
\newblock {Rigidity of {F}urstenberg entropy for semisimple {L}ie group actions},
\newblock {\em Ann. Sci. {\'E}cole Norm. Sup. (4)}, \textbf{33} (2000), 321--343.

\bibitem{NevoZ02}
\newblock A. Nevo and R.~J. Zimmer,
\newblock {A structure theorem for actions of semisimple {L}ie groups},
\newblock {\em Ann. of Math. (2)}, \textbf{156} (2002), 565--594.

\bibitem{Ontaneda:2020aa}
\newblock P. Ontaneda,
\newblock {Riemannian hyperbolization},
\newblock {\em Publ. Math. Inst. Hautes \'{E}tudes Sci.}, \textbf{131} (2020), 1--72.

\bibitem{MR2481335}
\newblock T.~L. Payne,
\newblock {Anosov automorphisms of nilpotent {L}ie algebras},
\newblock {\em J. Mod. Dyn.}, \textbf{3} (2009), 121--158.

\bibitem{MR1946555}
\newblock L. Polterovich,
\newblock {Growth of maps, distortion in groups and symplectic geometry},
\newblock {\em Invent. Math.}, \textbf{150} (2002), 655--686.

\bibitem{Prasad-Raghunathan1972}
\newblock G. Prasad and M.~S. Raghunathan,
\newblock {Cartan subgroups and lattices in semi-simple groups},
\newblock {\em Ann. of Math. (2)}, \textbf{96} (1972), 296--317.

\bibitem{Prasad-Rapinchuk2003}
\newblock G. Prasad and A.~S. Rapinchuk,
\newblock {Existence of irreducible {$\Bbb R$}-regular elements in {Z}ariski-dense subgroups},
\newblock {\em Math. Res. Lett.}, \textbf{10} (2003), 21--32.

\bibitem{MR2150880}
\newblock G. Prasad and A.~S. Rapinchuk,
\newblock {Zariski-dense subgroups and transcendental number theory},
\newblock {\em Math. Res. Lett.}, \textbf{12} (2005), 239--249.

\bibitem{raghunathan1972}
\newblock M.~S. Raghunathan,
\newblock \emph{Discrete Subgroups of Lie Groups}, Ergeb. Math. Grenzgeb., vol.~68,
\newblock Springer New York, 1972.

\bibitem{RH-W}
\newblock F. Rodriguez~Hertz and Z. Wang,
\newblock {Global rigidity of higher rank abelian {A}nosov algebraic actions},
\newblock {\em Invent. Math.}, \textbf{198} (2014), 165--209.

\bibitem{ShubExpandingMaps1970} 
\newblock M. Shub,
\newblock {Expanding maps},
\newblock in {\em Global {A}nalysis ({P}roc. {S}ympos. {P}ure {M}ath., {V}ol. {XIV}, {B}erkeley, {C}alif., 1968)},  273--276, Amer. Math. Soc., Providence, RI, 1970.

\bibitem{Smale}
\newblock S.~Smale,
\newblock {Differentiable dynamical systems},
\newblock {\em Bull. Amer. Math. Soc.}, \textbf{73} (1967), 747--817.

\bibitem{brinprize2015}
\newblock R. Spatzier,
\newblock{On the work of {R}odriguez {H}ertz on rigidity in dynamics},
\newblock {\em J. Mod. Dyn.}, \textbf{10} (2016), 191--207.

\bibitem{Spatzier:2020tx}
\newblock R.~Spatzier and K.~Vinhage,
\newblock {Cartan actions of higher rank abelian groups and their classification},
\newblock \emph{J. Amer. Math. Soc.}, \textbf{37} (2024),   731--859.

\bibitem{Spatzier-Yang2017}
\newblock R. Spatzier and L. Yang,
\newblock {Exponential mixing and smooth classification of commuting expanding maps},
\newblock {\em J. Mod. Dyn.}, \textbf{11} (2017), 263--312.

\bibitem{MR806415}
\newblock D. Sullivan,
\newblock {Quasiconformal homeomorphisms and dynamics. {II}. {S}tructural stability implies hyperbolicity for {K}leinian groups},
\newblock {\em Acta Math.}, \textbf{155} (1985), 243--260.

\bibitem{Vinhage:2022}
\newblock K. Vinhage,
\newblock {Instability for rank one factors of product actions},
\newblock J. Mod. Dyn., \textbf{21} (2025), 607--620, 2025.

\bibitem{Zimmer-Annals1980}
\newblock R.~J. Zimmer,
\newblock {Strong rigidity for ergodic actions of semisimple {L}ie groups},
\newblock {\em Ann. of Math. (2)}, \textbf{112} (1980), 511--529.

\bibitem{Zimmer-ICM}
\newblock R.~J. Zimmer,
\newblock Actions of semisimple groups and discrete subgroups,
\newblock in {\em Proceedings of the {I}nternational {C}ongress of {M}athematicians, {V}ol. {I, II} ({B}erkeley, {C}alif., 1986)}, 1247--1258, Amer. Math. Soc., Providence, RI, 1987.

\bibitem{Zimmer1}
\newblock R.~J. Zimmer,
\newblock {Lattices in semisimple groups and invariant geometric structures on compact manifolds},
\newblock in {\em Discrete Groups in Geometry and Analysis ({N}ew {H}aven, {C}onn., 1984)}, 152--210,  {Progr. Math.}, vol. 67, Birkh{\"a}user Boston, Boston, MA, 1987. 


\end{thebibliography}

\end{document}